\documentclass[11pt]{amsart}

\usepackage{amsfonts,amssymb,amscd,amsmath,enumerate,verbatim,calc}
\usepackage{amsmath}
\usepackage{amssymb}
\usepackage{amscd}

\usepackage{color}

\newcommand{\rmn}[1]{\romannumeral#1}
\def\RMN#1{\uppercase\expandafter{\romannumeral#1}}

\makeatletter
  
  \@addtoreset{equation}{section}
\makeatother

\newtheorem{Theorem}{Theorem}[section]
\newtheorem{Lemma}[Theorem]{Lemma}

\newtheorem{Proposition}[Theorem]{Proposition}
\newtheorem{Remark}[Theorem]{Remark}

\newtheorem{Example}[Theorem]{Example}
\newtheorem{Claim}[Theorem]{Claim}
\newtheorem{Definition}[Theorem]{Definition}

\newtheorem{Conjecture}[Theorem]{Conjecture}
\newtheorem{Question}[Theorem]{Question}
\newtheorem{Fact}[Theorem]{Fact}

\newcommand{\proj}[1]{\operatorname{Proj}(#1)}

\newcommand{\spec}{\operatorname{Spec}}

\newcommand{\rank}{\operatorname{rank}}

\newcommand{\A}{\operatorname{A}}

\newcommand{\ch}{\operatorname{ch}}

\theoremstyle{plain}

\title[The cone spanned by maximal Cohen-Macaulay modules]{The cone
  spanned by maximal Cohen-Macaulay modules and an application}
\author{C-Y. Jean Chan and Kazuhiko Kurano} 
\date{July~14, 2015. The final version, to appear in {\it
    Trans. Amer. Math. Soc.}.} 
\thanks{2010 Mathematics Subject Classifications: 13C14, 13D15, 13D40, 
14C17, 14C40.\\
\indent  Key words and phrases. Numerical rational equivalence,
  Cohen-Macaulay cone, test module, Hilbert-Kunz function, Segre product.\\
\indent The first author was partially supported by Early Career
  Investigator's Grant \#C61368 of Central Michigan University.  The
  second author was partially supported by KAKENHI (24540054). }

\begin{document} 

\begin{abstract}
  The aim of this paper is to define the notion of the
  {\em Cohen-Macaulay cone} of a Noetherian local domain $R$ and to present
 its applications to the theory of  Hilbert-Kunz functions. 
 It has been
  shown in \cite{K23} that, with a mild condition on  $R$, the 
  Grothendieck group $\overline{G_0(R)}$ 
of finitely generated $R$-modules modulo numerical equivalence 
is a finitely generated torsion-free abelian group.
 The Cohen-Macaulay cone of $R$ is the  cone
  in $\overline{G_0(R)}_{\mathbb R}$ spanned by cycles represented by 
maximal Cohen-Macaulay modules.
We study basic properties on the Cohen-Macaulay cone in this paper.
As an application,
  various examples of Hilbert-Kunz functions in the polynomial type
  will be produced. 
Precisely, for any given integers $\epsilon_i = 0, \pm 1$ ($d/2 < i < d$),
we shall construct a $d$-dimensional Cohen-Macaulay local ring $R$ 
(of characteristic $p$) and 
a maximal primary ideal $I$ of $R$ such that the function $\ell_R(R/I^{[p^n]})$
is a polynomial in $p^n$ of degree $d$ whose coefficient of $(p^n)^i$ is the
product of  $\epsilon_i$ and a positive rational number for $d/2 < i < d$. 
The existence of such ring is proved by using Segre
  products to construct a Cohen-Macaulay ring such that the Chow group of
  the ring is of certain simplicity and that {\em test modules} exists
  for it.
\end{abstract}

\maketitle

\section{Introduction}

Let $R$ be a Noetherian local domain.  In this paper,
we introduce the notion of the {\em Cohen-Macaulay cone} and 
{\em test modules} of $R$.  
As an application, for any given integers $\epsilon_i = 0, \pm 1$ ($d/2 < i < d$),
we shall construct a $d$-dimensional Cohen-Macaulay local ring $R$ 
(of characteristic $p$) and 
a maximal primary ideal $I$ of $R$ such that the function $\ell_R(R/I^{[p^n]})$
is a polynomial in $p^n$ of degree $d$ whose coefficient of $(p^n)^i$ is the
product of  $\epsilon_i$ and a positive rational number for $d/2 < i < d$. 

The main materials are divided into three parts. First part is an
introduction to the theory of the Cohen-Macaulay cone and test modules. 
Then we
prove the existence of a Cohen-Macaulay ring such that 
its Chow group of the ring is of certain simplicity
and it has a test module. 
Using the ring just constructed, we
shall further build a ring whose Hilbert-Kunz function satisfies 
the required conditions.

We now describe in more detail these new notions, and in the case
of positive characteristic, their contribution to the theory of
Hilbert-Kunz functions.

Let $R$ be a Noetherian local domain of dimension $d$. 
The Grothendieck group $\overline{G_0(R)}$ of 
finitely generated $R$-modules modulo numerical equivalence 
is defined and
studied in \cite{K23} where it is proven that under a mild
condition $\overline{G_0(R)}$ is a finitely generated torsion-free abelian group
(see also Theorem~\ref{K} in Section~2). 
Let
$\rho(R)$ denote the rank of
$\overline{G_0(R)}$.  In this paper, we let $R$ be a
Cohen-Macaulay local domain and introduce a cone inside 
${\mathbb R}^{\rho(R)} = \overline{G_0(R)} \otimes_{\mathbb Z} \mathbb R$
consisting of all nonnegative linear combinations of maximal Cohen-Macaulay
modules.  
This is called the {\em Cohen-Macaulay cone} of $R$. 
A module $M$
is a {\em test module} if $M$ is a maximal Cohen-Macaulay module such
that its Todd class consists
of only the top term; {\em i.e.}, $\overline{\tau_R}([M]) \in
\overline{{\rm A}_d(R)}_{\mathbb Q}$.  
In the case where $R$ is an F-finite Cohen-Macaulay local ring
(of positive characteristic $p$) with algebraically closed residue class field,
$M$ is a test module if and only if 
it is a maximal Cohen-Macaulay module such that
\[
[{}_{F^e}M] = p^{de}[M] \mbox{\ \ in $\overline{G_0(R)}_{\mathbb Q}$} 
\]
for some (any) $e > 0$,
where $[{}_{F^e}M]$ denotes the $R$-module $M$ whose $R$-module structure
is given by the $e$-th power of the Frobenius map.
If the small Cohen-Macaulay conjecture is affirmative,
then any Noetherian local ring has a test module.
However, the authors do not know whether 
test modules exist or not even if $R$ is a Gorenstein ring.
We refer the readers to \cite{K15} for test modules;
however, the definition of test modules in this paper is
slightly different from that in \cite{K15}.
In this paper, we need test modules which contain
the ring as a direct summand, that is,
test modules which are in the interior of the Cohen-Macaulay cone of $R$. 
The ideas of these new notions are motivated by
the studies of Hilbert-Kunz functions.

Assume additionally that $R$ has a positive characteristic $p$ with
dimension $d$. Let $I$ be a maximal primary ideal.  The Hilbert-Kunz
function with respect to $I$, named by Monsky~\cite{M}, is defined as
\[ \varphi_R(n)= \mbox{ length }({R/I^{[p^n]}R}) \] where $I^{[p^n]}$
is the Frobenius $n$-th power of $I$. Unlike the usual Hilbert
function, the shape of the Hilbert-Kunz function varies from case to
case.  Monsky proved that $\varphi_M(n) = e_{HK}(I,M) \, p^{nd} +
\mathcal O( p^{n(d-1)} )$ for some positive constant $ e_{HK}(I,M)$.
Its stability has been studied in Huneke-McDermott-Monsky~\cite{HMM},
Fakhruddin-Trivedi~\cite{FT}, Brenner~\cite{B},
Hochster-Yao~\cite{HY}, Chan-Kurano~\cite{CK1}, etc.

Classically from Macaulay's theorem one knows what numerical functions
are Hilbert functions ({\em c.f.} \cite[Section~4.2]{BH}). Similarly
for the Hilbert-Kunz function, it is natural to ask what functions are
Hilbert-Kunz functions.
 But the latter is a much more subtle question
since the shape of a Hilbert-Kunz function is unpredictable in general.
In \cite[Example 3.1(3)]{K24}, the second author proved that if
$I$ is a maximal primary ideal of a local ring $R$ 
that satisfies the following two conditions
 \begin{itemize} 
   \item $R$ is an $F$-finite Cohen-Macaulay local 
ring whose residue class field is algebraically  closed, and 
   \item $I$ has finite projective dimension,
  \end{itemize} 
then the Hilbert-Kunz function of $R$ with respect to $I$ is a polynomial of
  $p^n$ (see also \cite{K23} and \cite[Lemma 3.4]{CK1}).
We refer the reader to MacDonnell~\cite{MacDonnell} in the case where
$I$ is a homogeneous ideal.
One aim of this paper is to prove the following theorem:

\begin{Theorem}\label{HK}
Let $d$ be a positive integer and $p$ a prime number.
Let $\epsilon_0$, $\epsilon_1$, \ldots, $\epsilon_d$ be integers such that
\[
\epsilon_i = 
\left\{
\begin{array}{ll}
1 & i = d , \\
\mbox{$-1$, $0$ or $1$} & d/2 < i < d , \\
0 & i \le d/2 .
\end{array}
\right.
\]

Then, there exists a $d$-dimensional Cohen-Macaulay local ring $R$ of
characteristic $p$, a maximal primary ideal $I$ of $R$ of finite projective dimension,
and positive rational numbers $\alpha$, $\beta_{d-1}$, $\beta_{d-2}$,\ldots, $\beta_{0}$ such that
\[
\ell_R(R/I^{[p^n]}) = \epsilon_d \alpha p^{dn} + \sum_{i = 0}^{d-1} \epsilon_i \beta_i p^{in}
\]
for any $n > 0$.
\end{Theorem}

  The proof of Theorem~\ref{HK} is established by constructing a
  test module $M$ over a Cohen-Macaulay ring whose
  rational Chow group $\A_*(R)_{\mathbb Q}$ is of certain simplicity.
To determine how the Todd class of a module
  looks  is very difficult in general since $G_0(R)_{\mathbb Q}$ is too big. 
As mentioned earlier
  $\overline{G_0(R)}_{\mathbb R}$ is a finite dimensional 
 ${\mathbb R}$-vector space
  of dimension $\rho(R)$. 
We denote the Cohen-Macaulay cone  $\sum_{M: \mbox{\tiny
      MCM}} \mathbb R_{\geq 0} [M]$ in 
$\mathbb R^{\rho (R)} = \overline{G_0(R)}_{\mathbb R}$ by
$C_{CM}(R)$.
Then we prove that the existence of a test module is equivalent to certain properties 
on the projections of the Cohen-Macaulay cone $C_{CM}(R)$ in the Chow group
${\rm A}_*(R)_{\mathbb R}$ via the Riemann-Roch map
$\overline{\tau_R}$ (Theorem~\ref{main}). 

Assume that $R$ is a Cohen-Macaulay local domain.
If $\dim R \le 2$, then $\rho(R) = 1$.
If $\dim R \ge 3$, then there is no upper bound  for $\rho(R)$.
In either case, the Cohen-Macaulay cone has the
maximal possible dimension, {\em i.e.}, 
$\dim_{\mathbb R} C_{CM}(R) = \rho(R)$. 
Indeed if $\rho(R) = 1$, then $C_{CM}(R)$ is obviously
a half line. 
For arbitrary $\rho(R)$, we prove that there is an open
neighborhood $U$ of $[R]$ in $\mathbb R^{\rho (R)}$ such that $U$ is
contained in the interior of  $C_{CM}(R)$. 
This is proved in Lemma~\ref{Remark1.5} along with other
properties of the Cohen-Macaulay cone.

If the ring $R$ is of finite Cohen-Macaulay type---namely, there exist
only finitely many indecomposable isomorphism classes of maximal
Cohen-Macaulay modules---then the cone is finitely generated and so
it is closed in $\mathbb R^{\rho (R)}$ under the usual topology for
the Euclidean space. The Cohen-Macaulay cone in general may not be
finitely generated, but the authors do not know a Cohen-Macaulay ring
$R$ whose $C_{CM}(R)$ is not closed.

In order to know how the notion of Cohen-Macaulay cones and test modules
are applied in the study of Hilbert-Kunz functions, we provide a conceptual sketch
of a key step in the proof of Theorem~\ref{HK}. The idea presented below is 
loosely about  
Step 2.

By the singular Riemann-Roch theorem, the coefficients of the Hilbert-Kunz function with 
respect to $I$ of finite projective dimension over a Cohen-Macaulay ring can be 
expressed in terms of the localized Chern characters. 
Precisely it says that the coefficient of $(p^n)^i$ 
is $\ch_i(\mathbb G_{\bullet})(\tau_R([R])_i)$ where $\mathbb G_{\bullet}$ is the resolution
of $R/I$ and $\tau_R([R])_i$  in ${\rm A}_i(R)_{\mathbb R}$ is the $i$-th Todd class of $R$. 
Thus to obtain desired coefficients for the 
Hilbert-Kunz function is equivalent to obtaining the values of the corresponding  
localized Chern characters when applied to the Todd classes of $R$.  

We prove Theorem~\ref{HK} by constructing a module over a
Cohen-Macaulay ring $A$ such that the $i$-th Todd class of the module
in $\A_i(A)_{\mathbb Q}$ for each $i$ has the desired value when
the localized Chern character is applied to it. Then we take the idealization
of the module to obtain a ring whose Hilbert-Kunz function has the reqired form.
The module just described is constructed by induction on $i$, so the 
initial step that shows the existence of a ring $A$ that possesses
certain properties and a test module is crucial. This is done in
Lemma~\ref{lemma}.

Last, we would like to make a remark without the intension of getting into
any technical detail in the present paper. In Theorem~\ref{HK}, the coefficients of the
polynomial are assumed to be zero in the terms of degree $d/2$ or
lower. This assumption is made due to the fact that
$\overline{{\rm A}_i(R)}_{\mathbb Q} = 0$ for $i\leq d/2$ 
(for a homogeneous coordinate ring $R$ of a smooth projective variety) if the
Grothendieck's standard conjecture holds ({\em c.f.}
\cite[Remark~7.12]{K23}). There is no known example where $\overline{{\rm A}_i(R)}_{\mathbb Q}$ does
not vanish for some  $i \leq d/2$. 

The paper is organized as follows. Section~2 introduces the
Cohen-Macaulay cone of an arbitrary local domain and the definition of
test modules.  Basic properties of these new notions are proved.  In
Example~\ref{example}, assuming $R$ is complete or essentially of
finite type over a field, we construct some examples of
test modules.  We also prove equivalence conditions of the existence
of test modules; for a general local domain, it can be found in
Remark~\ref{rmkTest}. Main results about such equivalence for
Cohen-Macaulay rings are stated and proven in Theorem~\ref{main}.

Section~3 discusses the Hilbert-Kunz functions and proves
Theorem~\ref{HK}, which constructs rings whose Hilbert-Kunz function is
of the required form.  For Theorem~\ref{HK}, we need Lemma~\ref{lemma},
which assures the existence of a Cohen-Macaulay ring such that it has
test modules and its numerical Chow group satisfies certain
properties.  The proof of Lemma~\ref{lemma} deserves independent
attention and so Section~4 is devoted to proving this lemma. The Chow
group of $X=\mathbb P^m \times \mathbb P^n$ and the Riemann-Roch map
$\tau_X: G_0(X)_{\mathbb Q} \rightarrow {\rm A}_*(X)_{\mathbb Q}$ have
been carefully studied by the second
author~\cite{K11}. Lemma~\ref{lemma} is proven by taking appropriate
Segre products of graded rings and utilizing the special properties
on $\tau_X$ and ${\rm A}_*(X)_{\mathbb Q}$.


\bigskip

\section{Cohen-Macaulay cone}

Let $(R,m)$ be a $d$-dimensional Noetherian local domain.
We always assume that local domains are homomorphic images 
of regular local rings and assume that 
one of two conditions in Theorem~\ref{K} is satisfied.

Further, in this paper, we assume that 
all modules are finitely generated. 
Let $G_0(R)$ be the Grothendieck group of finitely generated $R$-modules.

We put
\[
C(R) = \left\{ {\Bbb F}. \ \left|
\begin{array}{l}
\mbox{bounded complex of finite $R$-free modules,} \\
\mbox{$H_i({\Bbb F}.)$ has finite length for any $i$}
\end{array}
\right. \right\} .
\]

For ${\Bbb F}. \in C(R)$, we define an additive map
\[
\chi_{{\Bbb F}.} : G_0(R) \longrightarrow {\Bbb Z}
\]
by 
\[
\chi_{{\Bbb F}.}([M]) = \sum_i(-1)^i \ell_R(H_i({\Bbb F}.\otimes M)) .
\]
We set
\[
\overline{G_0(R)} = 
G_0(R)/\{ c \in G_0(R) \mid 
\mbox{$\chi_{{\Bbb F}.}(c) = 0$ for any ${\Bbb F}. \in C(R)$} \}  .
\]

\begin{Theorem}[Kurano, \cite{K23} Theorem~3.1 and Remark~3.5]
\label{K}
Assume that a Noetherian local domain $R$ satisfies one of the 
following two conditions:
\begin{itemize}
\item
$R$ is an excellent ring such that $R$ contains ${\Bbb Q}$.
\item
$R$ is essentially of finite type over a field, ${\Bbb Z}$ or a
complete discrete valuation ring.
\end{itemize}
Then, $\overline{G_0(R)}$ is a finitely generated free ${\Bbb Z}$-module.
\end{Theorem}

If $\spec(R)$ has a resolution of singularities or a regular alteration, then the above theorem is still true for such $R$ without assuming one of two conditions above.

Let $\rho(R)$ be the rank of $\overline{G_0(R)}$.
Note that $\rho(R)$ is always positive, that is $\overline{G_0(R)}
\neq 0$. 
Indeed consider the Koszul complex $\Bbb K.$ of some system of
parameters $\underline a$ of $R$. 
Then $\chi_{\Bbb K.}([R])$ is equal to the Hilbert-Samuel multiplicity of
$R$ with respect to the ideal generated by $\underline{a}$ ({\em
  c.f.} \cite{Se} Chapter IV Theorem~1).
So $\chi_{\Bbb K.}([R]) \neq 0$. This shows that $[R]$ is not zero in
$\overline{G_0(R)}$ by definition ({\em c.f.} \cite{K23} page~582).
If $d \le 2$, then $\rho(R) = 1$ (see \cite{K23} Proposition~3.7).
For any given $d \ge 3$, there is no upper bound for $\rho(R)$ (see
\cite{K23} Example~4.1).

\begin{Proposition}
\begin{rm}
The following conditions are equivalent:
\begin{itemize}
\item[(1)]
$\rho(R) = 1$.
\item[(2)]
$\overline{G_0(R)} = {\Bbb Z}[R]$.
\item[(3)]
For any ${\Bbb F}. \in C(R)$ and any $R$-module $M$ with $\dim M < d$,
$\chi_{{\Bbb F}.}([M]) = 0$.
\item[(4)]
For any ${\Bbb F}. \in C(R)$ and any  
$R$-module $M$,
\[
\chi_{{\Bbb F}.}([M]) = {\rm rank}_RM \cdot \chi_{{\Bbb F}.}([R]) .
\]
\end{itemize}
\end{rm}
\end{Proposition}

\proof
It is easy to see $(4) \Longrightarrow (3) \Longrightarrow (2) \Longrightarrow (1)$.

We shall prove $(1) \Longrightarrow (4)$.  Let ${\Bbb K}.$ be the Koszul
complex of some system of parameters
$\underline a$. 
Then by
  Serre's theorem ({\em c.f.}  \cite{Se} Chapter IV Theorem~1),
  $\chi_{\Bbb K.}([M]) = e_{I}(M) $ where $e_I(-)$ denotes the
  Hilbert-Samuel multiplicity with respect to the ideal $I$ generated
  by $\underline a$. Since $e_I(M) = {\rm rank}_RM \cdot e_I(R)$, we
  have 
\[
\chi_{{\Bbb K}.}([M]) = {\rm rank}_RM \cdot \chi_{{\Bbb K}.}([R]) .
\]
On the other hand,
  note that $\chi_{\Bbb K.}([R]) =e_I(R) \neq 0$. Therefore, $[R]\neq
  0$ in $\overline{G_0(R)}_{\Bbb Q}$. Thus $\overline{G_0(R)}_{\Bbb Q}
  = {\Bbb Q}[R]$ by the condition (1). We write $[M]=r[R]$ in
  $\overline{G_0(R)}_{\Bbb Q}$ for some rational number $r$. Thus for
  every $\Bbb F. \in C(R)$, $\chi_{\Bbb F.}( [M]) = r \cdot \chi_{\Bbb
    F.}( [R] )$. In particular,
\[ r \cdot \chi_{\Bbb K.}( [R] ) = \chi_{\Bbb K.}( [M]) = {\rm rank}_R M
\cdot \chi_{\Bbb K.}( [R]) . \] This shows $r={\rm rank}_R M$ and $[M]
= {\rm rank}_RM \cdot [R]$ in $\overline{G_0(R)}_{\Bbb Q}$.
Therefore, the condition (4) is satisfied.  \qed

\begin{Remark}\label{1.3}
\begin{rm}
Suppose ${\Bbb F}. \in C(R)$.
Assume that ${\Bbb F}.$ is not exact 
and the length of ${\Bbb F}.$ is $d$, that is,
\[
{\Bbb F}. \ : \ 
0 \rightarrow F_d \rightarrow F_{d_1} \rightarrow 
\cdots \rightarrow F_1 \rightarrow F_0 \rightarrow 0 .
\]
Let $M$ be a finitely generated $R$-module. By induction on the depth,
we have $\operatorname{depth}( M) = d - \operatorname{max} \{ i | H_i(
{\Bbb F}. \otimes_R M ) \neq 0 \}$. 
If $M$ is a maximal Cohen-Macaulay module, then
$H_i({\Bbb F.}\otimes_RM) = 0$ for $i > 0$.  Thus,
\[
\chi_{{\Bbb F}.}([M]) = \ell(H_0({\Bbb F}.\otimes_RM)) > 0 .
\]
\end{rm}
\end{Remark}

We think that cycles in $\overline{G_0(R)}_{\Bbb R}$ 
represented by maximal Cohen-Macaulay modules are 
{\em positive} elements in a sense. 

In this paper, ${\mathbb Q}_{\ge 0}$ (resp.\ ${\mathbb R}_{\ge 0}$) denotes
the set of non-negative rational (resp.\ real) numbers.
Further, ${\mathbb Q}_+$ (resp.\ ${\mathbb R}_+$) denotes
the set of positive rational (resp.\ real) numbers.

\begin{Definition}
\begin{rm}
Set
\[
C_{CM}(R) = \sum_{M: MCM} {\Bbb R}_{\ge 0}[M] \subset 
\overline{G_0(R)}_{\Bbb R} ,
\]
where $M$ runs over all maximal Cohen-Macaulay $R$-modules
in the above summation.
We call it the {\em CM (Cohen-Macaulay) cone} of $R$.
Let $C_{CM}(R)^-$ be the closure of $C_{CM}(R)$
in $\overline{G_0(R)}_{\Bbb R}$ with respect to the usual topology on
the Euclidean space $\mathbb R^{\rho(R)}$.

The CM cone $C_{CM}(R)$ and its closure  $C_{CM}(R)^-$ are 
  convex cones by definition.  
The authors do not have an example where
$C_{CM}(R)$ is  a proper subset of its closure. 
(We know that $C_{CM}(R)^-$ is a strongly convex cone by a recent result
due to Dao and Kurano~\cite{DK2}.
We do not need this result in this paper.)

Set
\[
{\rm Nef}(R) = \{ c \in \overline{G_0(R)}_{\Bbb R} \mid
\mbox{$\chi_{{\Bbb F}.}(c) \ge 0$ for any ${\Bbb F}.
\in C(R)$ of length $d$}
\} .
\]
We call it the {\em nef (numerically effective) cone} of $R$.
\end{rm}
\end{Definition}

\begin{Lemma}\label{Remark1.5}
Let $R$ be a Cohen-Macaulay local domain.
\begin{enumerate}
\item
Let $c$ be in $\sum_{M: MCM}{\Bbb Q}_{\ge 0}[M]$.
Then, there exists a positive integer $n$ and a
maximal Cohen-Macaulay module $M$ such that
\[
nc = [M] \ \ \mbox{in $\overline{G_0(R)}_{\Bbb R}$.}
\]
\item
Let $c$ be an element in $C_{CM}(R)$.
For any open subset $U$ of $\overline{G_0(R)}_{\Bbb R}$
containing $c$,
\[
U
\cap \sum_{M: MCM}{\Bbb Q}_{\ge 0}[M] \neq \emptyset .
\]
\item
\[
C_{CM}(R) \cap \overline{G_0(R)}_{\Bbb Q} \subset \sum_{M: MCM}{\Bbb Q}_{\ge 0}[M] \subset \overline{G_0(R)}_{\Bbb R}
\]
is satisfied. 
\item
\[
C_{CM}(R) \subset C_{CM}(R)^- \subset {\rm Nef}(R) \subset 
\overline{G_0(R)}_{\Bbb R}.
\]
\item 
\[
{\rm Nef}(R) \cap - {\rm Nef}(R) = \{ 0 \} .
\]
\item
If $R$ is of finite Cohen-Macaulay representation type,
then
\[
C_{CM}(R) = C_{CM}(R)^- .
\]
\item
There exists an open set $U$ of $\overline{G_0(R)}_{\Bbb R}$ 
such that $[R] \in U \subset C_{CM}(R)$.
\item
\[
Int(C_{CM}(R)^-) \subset C_{CM}(R) .
\]
\end{enumerate}
\end{Lemma}

\proof
It is easy to see (1).

Here, we shall prove (2).
Put $c = \sum_ir_i[M_i]$, where $r_i \in {\Bbb R}_+$ and $M_i$ is
a maximal Cohen-Macaulay module for each $i$.
We choose $r'_i \in {\Bbb Q}_+$ sufficiently near $r_i$
for each $i$.
Then  $\sum_ir'_i[M_i]$ is in $U$.

We shall prove (3).
It is sufficient to show that, for a finite number of
maximal Cohen-Macaulay modules $M_1$, \ldots, $M_s$,
\[
\sum_{i = 1}^s{\Bbb R}_{\ge 0}[M_i]
\cap  \overline{G_0(R)}_{\Bbb Q}
\subset \sum_{i = 1}^s{\Bbb Q}_{\ge 0}[M_i]
\]
in $\overline{G_0(R)}_{\Bbb R}$.

Let $c = \sum_ih_i[M_i]$ ($h_i \in {\Bbb R}_+$) 
be in the left-hand side in the above.
First, remark that
\[
\sum_{i = 1}^s{\Bbb R}[M_i]
\cap  \overline{G_0(R)}_{\Bbb Q}
= \sum_{i = 1}^s{\Bbb Q}[M_i] .
\]
Therefore, 
\[
c = \sum_{i = 1}^sq_i[M_i]
\]
for $q_i \in {\Bbb Q}$ ($i = 1, \ldots, s$).
Here, we put
\[
W = \{ (\alpha_1, \ldots, \alpha_s) \in {\Bbb Q}^s \mid \sum_i\alpha_i[M_i] = 0 \}
\subset {\Bbb Q}^s .
\]
Then,
\[
(h_1 - q_1, \ldots, h_s - q_s)
\in W \otimes_{\Bbb Q}{\Bbb R}
\]
Therefore, there exists $(\beta_1, \ldots, \beta_s) \in W$
sufficiently near $(h_1 -q_1, \ldots, h_s - q_s)$.
Then,
\[
c = \sum_{i = 1}^s(q_i + \beta_i)[M_i]
\]
where $q_i + \beta_i \in {\Bbb Q}_+$ for $i = 1, \ldots, s$.

(4) immediately follows from Remark~\ref{1.3}.

We shall prove (5).
Let $c \in {\rm Nef}(R) \cap - {\rm Nef}(R)$.
For ${\Bbb F}. \in C(R)$ of length $d$, 
we have $\chi_{{\Bbb F}.}(c) \ge 0$ and $\chi_{{\Bbb F}.}(-c) \ge 0$.
Thus, we have $\chi_{{\Bbb F}.}(c) = 0$.
By Proposition~2 in \cite{RS},
the set of complexes of length $d$ generates the Grothendieck group of 
$C(R)$.
Therefore, $c$ is numerically equivalent to $0$.

(6) is easy.

Next, we shall prove (7).
Let $T_1$, \ldots, $T_\rho$ be torsion $R$-modules
such that
\begin{itemize}
\item
$\{ [T_1], \ldots, [T_{\rho - 1}], [R] \}$ is a basis
of the ${\Bbb Q}$-vector space $\overline{G_0(R)}_{\Bbb Q}$,
and
\item
$[T_1] + \cdots + [T_{\rho - 1}] + [T_\rho] = 0$
in $\overline{G_0(R)}_{\Bbb Q}$.
\end{itemize}
Let $M_i$ be the $k$-th sygyzy of $T_i$, where $k$ is an even integer 
bigger than $d$.
Then, $M_i$ is a maximal Cohen-Macaulay module such that
\[
[M_i] = ({\rm rank}_RM_i)[R] + [T_i]
\]
in $\overline{G_0(R)}_{\Bbb Q}$.
Then, we have
\[
[M_1] + \cdots + [M_{\rho - 1}] + [M_\rho] = 
\left( \sum_i{\rm rank}_RM_i \right) [R]
\]
in $\overline{G_0(R)}_{\Bbb Q}$.
By the above equations,
\[
\{ [M_1], \ldots, [M_{\rho - 1}],  [M_{\rho}] \}
\]
is a basis of the ${\mathbb R}$-vector space $\overline{G_0(R)}_{\Bbb R}$.
We know that $[R]$ is in interior of
the cone spanned by $[M_1]$, \ldots,
$[M_{\rho - 1}]$, $[M_\rho]$.

Lastly, we shall prove (8).
Suppose that $0 \neq c \in Int(C_{CM}(R)^-)$.
There exists an open neighborhood $U$ of $c$ 
such that $U$ is a subset of $C_{CM}(R)^-$.
Choose $e_1, \ldots, e_{\rho - 1} \in \overline{G_0(R)}_{\Bbb R}$ such that
$\{ e_1, \ldots, e_{\rho - 1}, c \}$ is an ${\Bbb R}$-basis of 
$\overline{G_0(R)}_{\Bbb R}$. 
Taking $e_i$'s small enough, we may assume $c + e_i \in U$ 
for $i = 1, \ldots, \rho - 1$, and $c - e_1 - \cdots - e_{\rho - 1} \in U$.
We put $s_i = c + e_i$ for $i = 1, \ldots, \rho - 1$ and
$s_\rho = c - e_1 - \cdots - e_{\rho - 1}$.
Then $s_1$, \ldots, $s_\rho$ are in $U$ such that
\[
s_1 +  \cdots + s_\rho = \rho \cdot c
\]
in  $\overline{G_0(R)}_{\Bbb R}$,
and $\{ s_1, \cdots, s_{\rho -1}, s_\rho \}$
is a basis of the ${\Bbb R}$-vector space 
$\overline{G_0(R)}_{\Bbb R}$.
Each point in $U$ is the limit of a sequence in 
$C_{CM}(R)$.
Therefore, there exist $s'_1, \ldots, s'_\rho$ such that
\begin{itemize}
\item
$s'_1, \ldots, s'_\rho \in C_{CM}(R)$,
\item
$\{ s'_1, \ldots, s'_\rho \}$ is a basis of the ${\mathbb R}$-vector space $\overline{G_0(R)}_{\Bbb R}$, and
\item
$c$ is in the cone spanned by $s'_1, \ldots, s'_\rho$.
\end{itemize}
Hence $c$ is in $C_{CM}(R)$.
\qed

\begin{Example}
\begin{rm}
Let
\[
R = k[x,y,z,w]_{(x,y,z,w)}/(xy-zw)
\]
Then,
\[
\{
R, P, Q \}
\]
is the set of isomorphism classes of indecomposable maximal Cohen-Macaulay modules (see Yoshino~\cite{Y}),
where $P = (x,z)$, $Q=(x,w)$.
In this case, $\rho(R) = 2$, and
\[
C_{CM}(R) = C_{CM}(R)^-
= {\mathbb R}_{\ge 0}[P] +  {\mathbb R}_{\ge 0}[Q]
 \subset {\rm Nef}(R) . 
\]
\end{rm}
\end{Example}

\begin{Remark}[Riemann-Roch theory]
\begin{rm}
We have an isomorphism of ${\Bbb Q}$-vector spaces
\[
\begin{array}{cccc}
G_0(R)_{\Bbb Q} & \stackrel{\tau_R}{\longrightarrow} & 
{\rm A}_*(R)_{\Bbb Q} & = \oplus_{i = 0}^d {\rm A}_i(R)_{\Bbb Q}
\end{array}
\]
as in \cite{F} and \cite{R}.
Then, we have ${\rm A}_d(R)_{\Bbb Q} = {\Bbb Q}[\spec(R)]$
and $p_d\tau_R([M]) = {\rm rank}_RM  \cdot [\spec(R)]$
where $p_d$ is the projection ${\rm A}_*(R) \rightarrow {\rm A}_d(R)$. 

Put $\tau_R([R]) = c_d + c_{d-1} + \cdots + c_0$,
where $c_i \in {\rm A}_i(R)_{\Bbb Q}$.
Then, 
\begin{enumerate}
\item
$c_d = [\spec(R)]$.
\item
If $R$ is a complete intersection,
$\tau_R([R]) = c_d$.
\item
If $R$ is Cohen-Macaulay,
\[
\tau_R([\omega_R]) = c_d - c_{d-1} + c_{d-2} - c_{d-3}+ \cdots ,
\]
where $\omega_R$ is the canonical module of $R$.
\item
If $R$ is Gorenstein,
\[
c_{d-1} = c_{d-3} = c_{d-5} = \ldots = 0 .
\]
\item
If $R$ is normal, we have an isomorphism ${\rm A}_{d-1}(R) \simeq {\rm Cl}(R)$ by
$[\spec(R/I)] \mapsto -{\rm cl}(I)$, where ${\rm cl}(I)$ denotes the isomorphism
class of a divisorial ideal $I$.
Then, we have
\[
c_{d-1} = - \frac{{\rm cl}({\omega_R})}{2} .
\]
\item
Localization of Galois extension of a regular local ring
satisfies $\tau_R([R]) = c_d$.
\end{enumerate}
\end{rm}
\end{Remark}

As in \cite{K23}, we can define $\overline{{\rm A}_i(R)}$
such that the following diagram is commutative:
\begin{equation}\label{zusiki}
\begin{array}{ccc}
{G_0(R)}_{\Bbb Q} & \stackrel{{\tau_R}}{\longrightarrow} & 
{{\rm A}_*(R)}_{\Bbb Q}  \\
\downarrow & & \downarrow  \\
\overline{G_0(R)}_{\Bbb Q} & \stackrel{\overline{\tau_R}}{\longrightarrow} & 
\overline{{\rm A}_*(R)}_{\Bbb Q}  \\
\downarrow & & \downarrow  \\
\overline{G_0(R)}_{\Bbb R} & \stackrel{\overline{\tau_R}}{\longrightarrow} & 
\overline{{\rm A}_*(R)}_{\Bbb R} 
\end{array}
\end{equation}

\begin{Definition}\label{test}
\begin{rm}
We say that $R$-module $M$ is an $R$-test module
if the following two conditions are satisfied:
\begin{enumerate}
\item
$M$ is a non-zero maximal Cohen-Macaulay module.
\item
$\overline{\tau_R}([M]) = \rank_RM \cdot [\spec(R)]$
in $\overline{{\rm A}_*(R)}_{\Bbb Q}$.
\end{enumerate}
\end{rm}
\end{Definition}

The above condition (2) is equivalent to $\overline{\tau_R}([M]) \in 
\overline{{\rm A}_d(R)}_{\Bbb Q}$.

The definition of test modules here is a little different from 
that in \cite{K15}.

For ${\Bbb F}. \in C(R)$, the Dutta multiplicity 
(limit multiplicity) is defined to be
\[
\chi_\infty ({\Bbb F}.) = \chi_{{\Bbb F}.}(\tau_R^{-1}([\spec(R)])) .
\]
If $M$ is an $R$-test module and ${\Bbb F}.$ is a complex in $C(R)$ of length $d$, then
\begin{equation}\label{positiveity}
\chi_\infty ({\Bbb F}.) = \frac{1}{{\rm rank}_RM}
\chi_{{\Bbb F}.}([M]) = \frac{1}{{\rm rank}_RM} \ell_R(H_0({\Bbb F}.\otimes_RM))
> 0 .
\end{equation}
For a non-exact complex in $C(R)$ of length $d$,
$\chi_{\infty}({\Bbb F}.)$ is positive if $R$ contains a field (\cite{Ri}, \cite{K13}, \cite{K15}).
Positivity of $\chi_{\infty}({\Bbb F}.)$ for a non-exact complex in $C(R)$
of length $d$ is an open question for the mixed characteristic case.
By (\ref{positiveity}), this question is true for $R$ which possesses a test module.

\begin{Remark}~\label{rmkTest} 
\begin{rm}
A local ring $R$ has a test module if and only if
\begin{equation}\label{defofTest}
\overline{\tau_R}^{-1}([\spec(R)]) \in C_{CM}(R)
\end{equation}
as follows.
The key point is that $\overline{\tau_R}^{-1}([\spec(R)])$
is in $\overline{G_0(R)}_{\Bbb Q}$ by the commutativity of the diagram~(\ref{zusiki}).
If (\ref{defofTest}) is satisfied, then we have
\[
\overline{\tau_R}^{-1}([\spec(R)]) \in 
C_{CM}(R) \cap \overline{G_0(R)}_{\Bbb Q} \subset \sum_{M: MCM}{\Bbb Q}_{\ge 0}[M]
\]
by Lemma~\ref{Remark1.5} (3).
Then, by Lemma~\ref{Remark1.5} (1),
we know the existence of $R$-test modules.
\end{rm}
\end{Remark}

\begin{Example}\label{example} 
\begin{rm}
Let $R$ be a Noetherian local domain of dimension $d$.
Suppose that $R$ contains an excellent regular local ring $S$,
and let $A$ be the integral closure of $S$ in $R$.
We assume that $A$ is a finitely generated $S$-module,
and $R$ coincides with $A_P$ for some prime ideal $P$ of $A$.
(We remark that such $S$ exists if $R$ is complete or essentially 
of finite type over a field.)

When the characteristic of $S$ is positive, we further assume that 
$S$ is F-finite.
Let $L$ be a finite dimensional normal extension of $Q(S)$
containing $Q(A)$ where $Q(S)$ and $Q(A)$ denote the field of fractions
of $S$ and $A$ respectively.
Let $B$ be the integral closure of $S$ in $L$.
Since $S$ is excellent, $B$ is a finitely generated $A$-module.
Thus $B\otimes_AR$ is a finitely generated $R$-module.
\begin{enumerate}
\item
Applying the method in \cite{K16}, we obtain
\[
\tau_R([B\otimes_AR]) \in {\rm A}_d(R)_{\Bbb Q} ,
\]
that is, 
\[
[B\otimes_AR] = {\rm rank}_R(B\otimes_AR) \cdot
\tau_R^{-1}([\spec(R)])
\]
in $G_0(R)_{\Bbb Q}$.
Therefore, if $B\otimes_AR$ is a Cohen-Macaulay ring,
then $B\otimes_AR$ is an $R$-test module.
\item
Put $G = {\rm Aut}_{Q(S)}(L)$.
Assume that $N$ is a maximal Cohen-Macaulay $B$-module.
For each $g \in G$, we give another $B$-module structure to $N$ by 
\[
\begin{array}{ccc}
B \times N & \longrightarrow & N \\
(b,n) & \mapsto & g(b)n . 
\end{array}
\]
We denote this $B$-module by ${}_gN$.
We put 
\[
M = \bigoplus_{g \in G}{}_gN .
\]
Then, we have 
\[
[M] = {\rm rank}_A(M) \cdot \tau_A^{-1}([\spec(A)])
\]
in $G_0(A)_{\Bbb Q}$.
Changing the base regular scheme using Lemma~4.1 (c) in \cite{K16},
we obtain
\[
[M\otimes_AR] = {\rm rank}_R(M\otimes_AR) \cdot 
\tau_R^{-1}([\spec(R)])
\]
in $G_0(R)_{\Bbb Q}$; therefore $M\otimes_AR$ 
is an $R$-test module.
\end{enumerate}
\end{rm}
\end{Example}

\begin{Remark}
\begin{rm}
\begin{enumerate}
\item
If any local ring has a test module,
then a conjecture (a positivity conjecture of Dutta multiplicity) 
is true (see (\ref{positiveity}),
\cite{K15} Conjecture~3.3 and Proposition~4.3). 
\item
Let $R$ be a complete local domain.
If the small Cohen-Macaulay conjecture is true,
then $R$ has a test module
 (see Example~\ref{example} and \cite{K15} Theorem~1.3). 
\item
Even if $R$ is a Gorenstein ring,
we do not know whether $R$ has a test module or not.
If $R$ is a complete intersection, then $R$ itself is an $R$-test module.
\end{enumerate}
\end{rm}
\end{Remark}

The following is the main theorem of this section.
We denote by $p_i$ the projection 
\[
\overline{{\rm A}_*(R)}_{\Bbb R} = 
\oplus_{i = 0}^d\overline{{\rm A}_i(R)}_{\Bbb R}
\rightarrow \overline{{\rm A}_i(R)}_{\Bbb R} .
\]

\begin{Theorem}\label{main}
Let $(R,m)$ be a Cohen-Macaulay local domain.
Consider the following conditions:
\begin{itemize}
\item[({\em \rmn{1}})]
$\overline{\tau_R}^{-1}([\spec(R)]) \in Int(C_{CM}(R)^-)$
\item[({\em \rmn{2}})]
$R$ has a test module which contains $R$ as a direct summand.
\item[({\em \rmn{3}})]
For $i = 0, 1, \ldots, d-1$, 
$p_i \overline{\tau_R}(C_{CM}(R)) = \overline{{\rm A}_i(R)}_{\Bbb R}$.
\item[({\em \rmn{4}})]
For $i = 0, 1, \ldots, d-1$, 
$p_i \overline{\tau_R}(C_{CM}(R)^-) = \overline{{\rm A}_i(R)}_{\Bbb R}$.
\end{itemize}

Then, we have
\[
(\rmn{1}) \Longleftrightarrow (\rmn{2}) \Longrightarrow (\rmn{3}) 
\Longleftrightarrow (\rmn{4}) .
\]

If $R$ is F-finite of characteristic $p > 0$ and $R/m$ is
algebraically closed, then above four conditions are equivalent
to each other.
\end{Theorem}

\proof
$(\rmn{1}) \Longrightarrow (\rmn{2})$.
There exists a positive integer $n$ such that
\[
n\overline{\tau_R}^{-1}([\spec(R)]) - [R] 
\in  Int(C_{CM}(R)^-) \cap \overline{G_0(R)}_{\Bbb Q}
\subset C_{CM}(R) \cap \overline{G_0(R)}_{\Bbb Q} 
\]
by Lemma~\ref{Remark1.5} (8).
By Lemma~\ref{Remark1.5} (1), (3), 
there exists a maximal Cohen-Macaulay 
module $M$ such that
\begin{equation}\label{s1}
[R] + [M] = n' \overline{\tau_R}^{-1}([\spec(R)]) 
\end{equation} 
for some $n' > 0$.

$(\rmn{2}) \Longrightarrow (\rmn{1})$.
Let $N_0$ be a module over $R$ such that $N = R\oplus N_0$ is a test
module. Let $M = N_0 \oplus N$. Then $M$ is a
maximal Cohen-Macaulay module, and we have an equality as (\ref{s1})
in which $n' = 2 \rank_R N$. 
Since $[R] \in Int(C_{CM}(R)^-)$ by Lemma~\ref{Remark1.5} (7) and
$[M] \in C_{CM}(R)$,  $[R \oplus M]$ is also in $Int(C_{CM}(R)^-)$.

$(\rmn{3}) \Longrightarrow (\rmn{4})$ is trivial.

$(\rmn{4}) \Longrightarrow (\rmn{3})$.
Since
\[
p_i\overline{\tau_R}(C_{CM}(R)^-) \subset (p_i\overline{\tau_R}(C_{CM}(R)))^- ,
\]
we have $(p_i\overline{\tau_R}(C_{CM}(R)))^- = \overline{{\rm
    A}_i(R)}_{\Bbb R}$.  Note that $p_i \overline{\tau_R}(C_{CM}(R))$
is a cone in $\overline{{\rm A}_i(R)}_{\Bbb R}$.  Since the convexity
is preserved under $\overline{\tau_R}$ and the projection $p_i$, 
if $p_i
\overline{\tau_R}(C_{CM}(R)) \neq \overline{{\rm A}_i(R)}_{\Bbb R}$,
then it must be contained in a closed half-space, and so must its
closure. This contradicts the above fact
$(p_i\overline{\tau_R}(C_{CM}(R)))^- = \overline{{\rm A}_i(R)}_{\Bbb
  R}$ resulted from the condition $(\rmn{4})$.  Hence
$p_i\overline{\tau_R}(C_{CM}(R)) = \overline{{\rm A}_i(R)}_{\Bbb R}$.

$(\rmn{1}) \Longrightarrow (\rmn{3})$.
Remark that
\[
p_i\overline{\tau_R}(C_{CM}(R)) \supset 
p_i\overline{\tau_R}(Int(C_{CM}(R)^-))
\ni p_i([\spec(R)]) = {\bf 0} 
\]
if $i < d$.
Since $p_i\overline{\tau_R}$ is an open map,
$p_i\overline{\tau_R}(Int(C_{CM}(R)^-))$ contains an open neighbourhood of ${\bf 0}$.
Then, we have $p_i\overline{\tau_R}(C_{CM}(R)) = \overline{{\rm A}_i(R)}_{\Bbb R}$
since $p_i \overline{\tau_R}(C_{CM}(R))$ is a cone in $\overline{{\rm A}_i(R)}_{\Bbb R}$.

Now, we shall prove 
$(\rmn{3}) \Longrightarrow (\rmn{1})$
in the case where $R$ is F-finite of characteristic $p > 0$ 
and $R/m$ is algebraically closed.

\noindent
{\bf Step 1} \
First we want to show $\overline{\tau_R}^{-1}([\spec(R)]) \in C_{CM}(R)^-$.
We put $\overline{\tau_R}([R]) = c_d + c_{d-1} + \cdots + c_0$,
where $c_i \in \overline{{\rm A}_i(R)}_{\Bbb Q}$.
\[
\begin{array}{ccl}
\overline{G_0(R)}_{\Bbb R} & \stackrel{\overline{\tau_R}}{\longrightarrow} &
\overline{{\rm A}_*(R)}_{\Bbb R} = \oplus_{i = 0}^d\overline{{\rm A}_i(R)}_{\Bbb R} \\
\mbox{$[R]$}  & \longleftrightarrow & c_d + c_{d-1} + \cdots + c_0 \\
\mbox{$[R^{\frac{1}{p^e}}]$} &\longleftrightarrow 
& p^{de}c_d + p^{(d-1)e}c_{d-1} + \cdots + p^{0e}c_0 .
\end{array}
\]
Therefore
\[
\frac{1}{p^{de}}[R^{\frac{1}{p^e}}] =
\overline{\tau_R}^{-1}\left( c_d + \frac{1}{p^e}c_{d-1} + \cdots + \frac{1}{p^{de}}c_0 \right)
\in C_{CM}(R) .
\]
Take the limit. Then, we have
\[
\overline{\tau_R}^{-1}([\spec(R)]) = \overline{\tau_R}^{-1}(c_d) = 
\lim_{e \rightarrow \infty} 
\frac{1}{p^{de}}[R^{\frac{1}{p^e}}]
\in C_{CM}(R)^- .
\]

\noindent
{\bf Step 2} \
We shall show $\overline{\tau_R}^{-1}([\spec(R)]) \in Int(C_{CM}(R)^-)$.

Assume that $\overline{\tau_R}^{-1}([\spec(R)])$ is in the boundary of
the cone $C_{CM}(R)^-$.  Then $[\spec R]$ is in the
  boundary of the image of the cone under $\overline{\tau_R}$.  If $\rho(R)=1$,
then it can never happen. Therefore we may assume that $\rho(R) >1$.
For any $R$-module $M$,
\[
\overline{\tau_R}([M]) = \rank_RM [\spec(R)]
+ (\mbox{lower dimensional terms}) .
\]
Therefore, $\overline{\tau_R}(C_{CR}(R)^-) \neq \overline{{\rm A}_*(R)}_{\Bbb R}$.
Note that
  $\overline{\tau_R} ( C_{CM}(R)^-)$ is a convex cone since
  $\overline{\tau_R}$ is an ${\Bbb R}$-linear map. Thus there exists a
    hyperplane through the origin that contains the boundary
    possessing $[\spec(R)]$. Indeed such a hyperplane is a supporting
    hyperplane of  the cone $\overline{\tau_R} (C_{CM}(R)^-)$; namely,
    there exists a vector ${\bf v}$ normal to the
    hyperplane such that the inner product $<\bf v, \bf u>$ is
    nonnegative for every $\bf u$ in $\overline{ \tau_R}(
    C_{CR}(R)^-)$. Let $\xi$ be the projection of $ \overline{{\rm
        A}_*(R)}_{\Bbb R}$ onto the line generated by $\bf v$. Then
\[
\xi : \overline{{\rm A}_*(R)}_{\Bbb R} \longrightarrow {\Bbb R}
\]
is a non-zero ${\Bbb R}$-linear map
with the properties
\begin{equation}\label{katei}
\left\{
\begin{array}{l}
\xi([\spec(R)]) = 0 , \\
\xi\overline{\tau_R}(C_{CM}(R)^-) \subset {\Bbb R}_{\ge 0} .
\end{array}
\right.
\end{equation}

Since $\overline{{\rm A}_d(R)}_{\Bbb R} = {\Bbb R}[\spec(R)]$,
we have $\xi(\overline{{\rm A}_d(R)}_{\Bbb R}) = 0$.
Since $\xi \neq 0$, we can choose $0 \le j < d$
such that
\begin{equation}\label{xi}
\left\{
\begin{array}{l}
\mbox{$\xi(\overline{{\rm A}_i(R)}_{\Bbb R}) = 0$ for $i = j+1, j+2, \ldots, d$}, \\
\xi(\overline{{\rm A}_{j}(R)}_{\Bbb R}) \neq 0 .
\end{array}
\right.
\end{equation}
Therefore, $\xi(\overline{{\rm A}_{j}(R)}_{\Bbb R}) = {\Bbb R}$.
Since $p_j\overline{\tau_R}(C_{CM}(R)) = \overline{{\rm A}_j(R)}_{\Bbb R}$
by the condition~(\rmn{3}), we have
\[
{\Bbb R} = \xi p_j\overline{\tau_R}(C_{CM}(R)) = \sum_{M: MCM} {\Bbb R}_{\ge 0}\xi p_j\overline{\tau_R}([M]) .
\]
Therefore, 
there exists a maximal Cohen-Macaulay module $N$ such that 
\begin{equation}\label{xip_jtau}
\xi p_j\overline{\tau_R}([N]) < 0 .
\end{equation}
Set
\[
\overline{\tau_R}([N]) = s_d + s_{d-1} + \cdots + s_0 ,
\]
where $s_i \in \overline{{\rm A}_i(R)}_{\Bbb Q}$.
By (\ref{xip_jtau}), we have
\begin{equation}\label{<}
\xi(s_j) < 0 .
\end{equation}
Then,
\[
\overline{\tau_R}([F^e(N)]) = 
p^{de}s_d + p^{(d-1)e}s_{d-1} + \cdots + p^{0e}s_0 .
\]
By the assumption~(\ref{xi}),
we have
\[
\xi\overline{\tau_R}([F^e(N)]) = 
p^{je}\xi(s_j) + p^{(j-1)e}\xi(s_{j-1}) + \cdots +
p^{0e}\xi(s_0) .
\]
Since (\ref{<}),
we know
\[
\xi\overline{\tau_R}([F^e(N)]) < 0
\]
for a sufficiently large $e$.
Since $F^e(N)$ is Cohen-Macaulay,
\[
\xi\overline{\tau_R}([F^e(N)]) \in \xi\overline{\tau_R}(C_{CM}(R))
\subset {\Bbb R}_{\ge 0} 
\]
by (\ref{katei}).
It is a contradiction.
\qed

For a positive integer $\ell$, we define
\[
\psi^\ell : \overline{{\rm A}_*(R)}_{\Bbb R} \longrightarrow \overline{{\rm A}_*(R)}_{\Bbb R}
\]
to be
\[
\psi^\ell(s_d + s_{d-1} + \cdots + s_0)
= \ell^ds_d + \ell^{d-1}s_{d-1} + \cdots + \ell^0s_0 ,
\]
where $s_i \in \overline{{\rm A}_i(R)}_{\Bbb R}$
for $i = 0, 1, \ldots, d$.

If there exists $\ell \ge 2$ such that
\[
\psi^\ell
\left(
\overline{\tau_R}(C_{CM}(R))
\right)
\subset
\overline{\tau_R}(C_{CM}(R)) ,
\]
the conditions (\rmn{1}), (\rmn{2}), (\rmn{3}), (\rmn{4})
in Theorem~\ref{main}
are equivalent to each other 
without assuming that $R$ is of positive characteristic.

If $R$ is of characteristic prime $p$,
then 
\[
\psi^p
\left(
\overline{\tau_R}(C_{CM}(R))
\right)
\subset
\overline{\tau_R}(C_{CM}(R)) .
\]

Therefore, it is natural to ask the following 
for an arbitrary Cohen-Macaulay local domain:

\begin{Question}
\begin{rm}
Is there an integer $\ell \ge 2$ such that
\[
\psi^\ell
\left(
\overline{\tau_R}(C_{CM}(R))
\right)
\subset
\overline{\tau_R}(C_{CM}(R)) ?
\]
\end{rm}
\end{Question}


\section{Examples of various Hilbert-Kunz functions}
\label{section2}
In the rest of this paper, we shall prove Theorem~\ref{HK} in the introduction.

We need the following lemma.

\begin{Lemma}\label{lemma}
Let $d$ be a positive integer, and $p$ be a prime number.

Then, there exists a ring $A$ 
of characteristic $p$ which satisfies the following conditions:
\begin{itemize}
\item
$A$ is a $d$-dimensional $F$-finite Cohen-Macaulay normal local domain and the residue class field of $A$ is algebraically closed.
\item
${\rm A}_i(A)_{\Bbb Q} = 
\overline{{\rm A}_i(A)}_{\Bbb Q} = 
\left\{
\begin{array}{lll}
{\Bbb Q} & & (\frac{d}{2} < i \le d) \\
0 & & (\mbox{otherwise})
\end{array}
\right.$.
\item
There exists a maximal Cohen-Macaulay $A$-module $M$ such that 
$\tau_A([A \oplus M]) \in {\rm A}_d(A)_{\Bbb Q}$; that is,
$A \oplus M$ is an $A$-test module containing $A$ as a direct summand.
\end{itemize}
\end{Lemma}

The above lemma will be proven in the next section.
In this section, using Lemma~\ref{lemma},
we shall prove Theorem~\ref{HK}.

Let $A$ be a ring satisfying three conditions in Lemma~\ref{lemma}.

\vspace{2mm}

\noindent
{\bf Step~1}.
We set
\[
\{ i_1, i_2, \ldots, i_t \} = \{ i \mid \epsilon_i \neq 0 \} .
\]
In Step~1, we shall show that
there exists $a_{i_k} \in {\rm A}_{i_k}(A)_{\Bbb Q}$ 
for $k = 1, 2, \ldots, t$, and a finite free $A$-complex ${\Bbb F}.$ of length $d$ with support at 
the maximal ideal $m$ such that
\[
{\rm ch}({\Bbb F}.) (a_{i_k}) \neq 0
\]
for all $k = 1, 2, \ldots, t$.

Recall that $\epsilon_j = 0$ if $j \le \frac{d}{2}$.
Then, by the assumption on the ring $A$, 
$\overline{{\rm A}_{i_k}(A)}_{\Bbb Q} = {\Bbb Q}$ for 
$k = 1, 2, \ldots, t$.
By the definition of $\overline{{\rm A}_{i}(A)}_{\Bbb Q}$ 
(see \cite{K23}), 
there exists $a_{i_k} \in {\rm A}_{i_k}(A)_{\Bbb Q}$ 
 and a finite free $A$-complex ${\Bbb F}^{(k)}.$ of length $d$ with support at 
the maximal ideal $m$ such that
\[
{\rm ch}({\Bbb F}^{(k)}.) (a_{i_k}) \neq 0
\]
for $k = 1, 2, \ldots, t$.
Here, we recall that, 
since $A$ is a Cohen-Macaulay local ring, 
the Grothendieck group of bounded finite $A$-free complexes 
with support in $\{ m \}$ is generated by free resolutions
of modules of finite length and of finite projective dimesion
(cf.\ Proposition~2 in \cite{RS}).

By induction, it is easy to show that there exist 
positive integers
$n_1$, $n_2$, \ldots, $n_t$ such that
\[
{\Bbb F}. = {{\Bbb F}^{(1)}.}^{\oplus n_1} \oplus {{\Bbb F}^{(2)}.}^{\oplus n_2} \oplus
\cdots \oplus {{\Bbb F}^{(t)}.}^{\oplus n_t}
\]
satisfies the required condition.

\vspace{2mm}

\noindent
{\bf Step~2}.
Take the complex ${\Bbb F}.$ that we have constructed in Step~1.

In Step~2, we shall show that
there exists a maximal Cohen-Macaulay $A$-module $N$ and positive rational numbers 
$\beta_0$, $\beta_1$, \ldots, $\beta_d$ such that
\[
{\rm ch}({\Bbb F}.) (\tau_A([A \oplus N])_i) = \epsilon_i \beta_i
\]
for $i = 0, 1, \ldots, d$.
Here, $\tau_A([A \oplus N])_i$ is the element in ${\rm A}_i(A)_{\Bbb Q}$ 
such that $\tau_A([A \oplus N]) = \sum_{i = 0}^d\tau_A([A \oplus N])_i$.

By the induction on $j$, we shall prove the following:
There exists a maximal Cohen-Macaulay $A$-module $N$ and positive rational numbers 
$\beta_{d-j}$, $\beta_{d-j+1}$, \ldots, $\beta_d$ such that
\[
{\rm ch}({\Bbb F}.) (\tau_A([A \oplus N])_i) = \epsilon_i \beta_i
\]
for $i = d-j, d-j+1, \ldots, d$.

Consider the case $j = 0$.
Here, recall that ${\Bbb F}.$ is a bounded finite $A$-free complex 
of length $d$ with support in $\{ m \}$.
Set $N = A$.
Then, $\tau_A([A \oplus A])_d = 2[{\rm Spec}(A)]$ and
\[
{\rm ch}({\Bbb F}.) ([{\rm Spec}(A)]) > 0
\]
by a theorem of Roberts \cite{Ri}. 
Here, recall that $\epsilon_d = 1$.
Therefore, $N=A$ 
and $\beta_d = {\rm ch}({\Bbb F}.) (2 [{\rm Spec}(A)])$ 
satisfy the required condition.

Next suppose $0 \le j < d$.
We assume that there exists a maximal Cohen-Macaulay $A$-module $N'$ 
and positive rational numbers 
$\beta'_{d-j}$, $\beta'_{d-j+1}$, \ldots, $\beta'_d$ such that
\[
{\rm ch}({\Bbb F}.) (\tau_A([A \oplus N'])_i) = \epsilon_i \beta'_i
\]
for $i = d-j, d-j+1, \ldots, d$.

Compare the rational number ${\rm ch}({\Bbb F}.) (\tau_A([A \oplus N'])_{d-j-1})$ with $\epsilon_{d-j-1}$.

If there exists a positive rational number $\beta'_{d-j-1}$ such that 
\[
{\rm ch}({\Bbb F}.) (\tau_A([A \oplus N'])_{d-j-1}) = \epsilon_{d-j-1} \beta'_{d-j-1} ,
\]
we have nothing to prove.
(Here, if $d - j - 1 \le \frac{d}{2}$, both ${\rm ch}({\Bbb F}.) (\tau_A([A \oplus N'])_{d-j-1})$ and $\epsilon_{d-j-1}$ are $0$.
Therefore, we have only to set $\beta'_{d-j-1} = 1$ in this case.)

We assume that there does not exist a positive rational number  $\beta'_{d-j-1}$
satisfying the above condition.
\begin{itemize}
\item[(*)]
If $\epsilon_{d-j-1} = 0$, we set $b = - \tau_A([A \oplus N'])_{d-j-1}$.
If $\epsilon_{d-j-1} \neq 0$, we choose $b \in {\rm A}_{d-j-1}(A)_{\Bbb Q}$
such that the sign of ${\rm ch}({\Bbb F}.)(b)$ is the same as that of $\epsilon_{d-j-1}$.
\end{itemize}
Here, remark that, by the construction of ${\Bbb F}.$ in Step~1,
we can choose an element $b$ satisfying the above condition.
We shall show the following claim:

\begin{Claim}\label{tyutyou}
There exists a maximal Cohen-Macaulay $A$-module $L$ and a positive integer $n$ 
such that 
\[
\tau_A([L]) = {\rm rank}_A(L) [{\rm Spec}(A)] + nb + (\mbox{lower dimensional terms}) .
\]
\end{Claim}

We prove this claim.

Since $b \in {\rm A}_{d-j-1}(A)_{\Bbb Q}$ and $d-j-1 < d$,
there exist (not necessary distinct) prime ideals $P_1$, \ldots, $P_s$
of height $j+1$ such that
\[
nb = [{\rm Spec}(A/P_1)] + \cdots + [{\rm Spec}(A/P_s)]
\]
in ${\rm A}_{d-j-1}(A)_{\Bbb Q}$ for some positive integer $n$.

Consider the following exact sequence
\[
0 \rightarrow N_1 \rightarrow F_{2 u-1} \rightarrow \cdots \rightarrow F_1 \rightarrow F_0 
\rightarrow A/P_1 \oplus \cdots \oplus A/P_s \rightarrow 0,
\]
where $F_0$, $F_1$, \ldots, $F_{2 u-1}$ are finitely generated $A$-free modules
and $u$ is a large enough number such that $N_1$ is a maximal Cohen-Macaulay $A$-module.

Then, we have
\begin{eqnarray*}
[N_1] & = & [A/P_1 \oplus \cdots \oplus A/P_s] - [F_0] + [F_1] - \cdots + [F_{2u-1}] \\
& = & [A/P_1 \oplus \cdots \oplus A/P_s] + {\rm rank}_A(N_1) [A] 
\end{eqnarray*}
in ${\rm G}_0(A)_{\Bbb Q}$.

By the assumption, there exists a maximal Cohen-Macaulay $A$-module $M$ such that 
$\tau_A([A \oplus M]) \in {\rm A}_d(A)_{\Bbb Q}$.
Adding ${\rm rank}_A(N_1) [M]$ to the both sides, we obtain
\[
[N_1] + {\rm rank}_A(N_1) [M] =  [A/P_1 \oplus \cdots \oplus A/P_s] + 
{\rm rank}_A(N_1) [A \oplus M] 
\]
in ${\rm G}_0(A)_{\Bbb Q}$.

Then, we have 
\begin{eqnarray*}
& & \tau_A([N_1 \oplus M^{\oplus {\rm rank}_A(N_1)}]) \\ & = &
{\rm rank}_A(N_1 \oplus M^{\oplus {\rm rank}_A(N_1)}) [{\rm Spec}(A)] + 
nb + (\mbox{lower dimensional terms})
\end{eqnarray*}
since 
\[
\tau_A( [A/P_1 \oplus \cdots \oplus A/P_s] )
= nb + (\mbox{lower dimensional terms})
\]
by the top term property (Theorem~18.3 (5) in \cite{F}).
Thus, $N_1 \oplus M^{\oplus {\rm rank}_A(N_1)}$ satisfies the condition on $L$ in Claim~\ref{tyutyou}.
We have completed the proof of Claim~\ref{tyutyou}.

\vspace{2mm}

Here, we set 
\[
N = A^{\oplus e - 1} \oplus N'^{\oplus e} \oplus 
 L^{\oplus f} .
\]
for some positive integers $e$ and $f$, and a maximal Cohen-Macaulay 
module $L$ in Claim~\ref{tyutyou}.
Then, 
\[
\tau_A([A \oplus N])_d = {\rm rank}_A(A \oplus N) [{\rm Spec}(A)]
= \frac{{\rm rank}_A(A \oplus N)}{{\rm rank}_A(A \oplus N')}
\tau_A([A \oplus N'])_d
\]
and
\[
\tau_A([A \oplus N])_i = e \tau_A([A \oplus N'])_i
\]
for $i = d-j, d-j+1, \ldots, d-1$.
Therefore, we have 
\[
{\rm ch}({\Bbb F}.)(\tau_A([A \oplus N])_d)
= \epsilon_d\beta_d
\]
where
\[
\beta_d = \frac{{\rm rank}_A(A \oplus N)}{{\rm rank}_A(A \oplus N')}
\beta'_d > 0 ,
\]
and
\[
{\rm ch}({\Bbb F}.)(\tau_A([A \oplus N])_i)
= \epsilon_i\beta_i
\]
where
\[
\beta_i = e \beta'_i > 0 
\]
for $i = d-j, d-j+1, \ldots, d-1$.

On the other hand, we have
\[
\tau_A([A \oplus N])_{d-j-1} = 
e \tau_A([A \oplus N'])_{d-j-1} + fnb .
\]

If $\epsilon_{d-j-1} = 0$, then we suppose $e = fn$. 
Then, $\tau_A([A \oplus N])_{d-j-1} = 0$ by the definition of $b$
(see (*) above Claim~\ref{tyutyou}).
Thus, putting $\beta_{d-j-1} = 1$,
\[
{\rm ch}({\Bbb F}.)(\tau_A([A \oplus N])_{d-j-1})
= 0 = \epsilon_{d-j-1}\beta_{d-j-1} .
\]

Next, assume that $\epsilon_{d-j-1} \neq 0$.
Consider the equality
\[
{\rm ch}({\Bbb F}.)(\tau_A([A \oplus N])_{d-j-1})
= e \ {\rm ch}({\Bbb F}.)(\tau_A([A \oplus N'])_{d-j-1})
+ fn \ {\rm ch}({\Bbb F}.)(b) .
\]
Assume that  $f/e$ is big enough.
Then the sign of the right-hand-side
is the same as that of $\epsilon_{d-j-1}$ by the definition of $b$
(see (*) above Claim~\ref{tyutyou}).
Therefore, there exists a positive rational number
$\beta_{d-j-1}$ such that
\[
{\rm ch}({\Bbb F}.)(\tau_A([A \oplus N])_{d-j-1})
= \epsilon_{d-j-1}\beta_{d-j-1} .
\]

\vspace{2mm}

\noindent
{\bf Step~3}.
Let ${\Bbb F}.$ and $N$ be a complex and 
a maximal Cohen-Macaulay module as in Step~2, respectively.
Let $R$ be the idealization $A \ltimes N$.
Then, $R$ is a $d$-dimensional Cohen-Macaulay local ring.

Since ${\rm ch}({\Bbb F}.)$ is a bivariant class (Chapter~17 in \cite{F}),
we have the commutative diagram
\[
\begin{array}{ccc}
{\rm A}_*(R/(m \ltimes mN))_{\Bbb Q} & 
\stackrel{{\rm ch}({\Bbb F}.\otimes_AR)}{\longleftarrow} &
{\rm A}_*(R)_{\Bbb Q} \\
\downarrow & & \downarrow \\
{\rm A}_*(A/m)_{\Bbb Q} & 
\stackrel{{\rm ch}({\Bbb F}.)}{\longleftarrow} &
{\rm A}_*(A)_{\Bbb Q}
\end{array}
\]
where the vertical maps are isomorphisms
induced by finite morphisms ${\rm Spec}(R) \rightarrow {\rm Spec}(A)$
and ${\rm Spec}(R/(m \ltimes mN)) \rightarrow {\rm Spec}(A/m)$.
Then, we obtain
\[
{\rm ch}({\Bbb F}.\otimes_AR)(\tau_R([R])_i) =
{\rm ch}({\Bbb F}.)(\tau_A([A \oplus N])_i) =
\epsilon_i \beta_i 
\]
for $i = 0, 1, \ldots, d$.
Since $R$ is a Cohen-Macaulay local ring of dimension $d$,
${\Bbb F}.\otimes_AR$ is a finite free resolution of an $R$-module 
$Q$ of finite length by \cite{RS}.
Let ${\mathcal C}$ be the category of $R$-modules of finite length 
and finite projective dimension.

Then, by Kumar's method (cf.\ Lemma 9.10 in \cite{Sr}),
there exist maximal primary ideals $I_1$, \ldots, $I_\ell$, $I$ of $R$
of finite projective dimension such that
\begin{itemize}
\item
$I_i$ is an ideal generated by a maximal regular sequence of $R$
for $i = 1, \ldots, \ell$.
\item
$[Q] + \sum_{i = 1}^\ell [R/I_i] = [R/I]$ in 
${\rm K}_0({\mathcal C})$.
\end{itemize}

Let $F : R \rightarrow R$ be the Frobenius map.
It is a finite morphism since $A$ is $F$-finite.
We denote by ${}_{F^n}R$ the $R$-module $R$ whose $R$-module structure
is given by 
\[
r \times a \mathrel{\mathop{:=}} r^{p^n}a 
\]
for $r \in R$ and $a \in {}_{F^n}R$.

By the Riemann-Roch formula,
\[
\tau_R([{}_{F^n}R]) = \sum_{i = 0}^dp^{in}\tau_R([R])_i .
\]
By the local Riemann-Roch formula, we have
\begin{eqnarray*}
\chi_{{\Bbb F}.\otimes_AR}({}_{F^n}R) 
& = & {\rm ch}({\Bbb F}.\otimes_AR)(\sum_{i = 0}^d\tau_R([{}_{F^n}R])_i) \\
& = & {\rm ch}({\Bbb F}.\otimes_AR)(\sum_{i = 0}^dp^{in}\tau_R([R])_i) \\
& = & \sum_{i = 0}^d {\rm ch}({\Bbb F}.\otimes_AR)(\tau_R([R])_i) p^{in} \\
& = & \sum_{i = 0}^d \epsilon_i \beta_i p^{in} .
\end{eqnarray*}

On the other hand, we have
\begin{eqnarray*}
\chi_{{\Bbb F}.\otimes_AR}({}_{F^n}R) & = & 
\chi(Q, {}_{F^n}R) \\
& = & \chi(R/I, {}_{F^n}R) - \sum_{i = 1}^\ell \chi(R/I_i, {}_{F^n}R) \\
& = & \ell_{R}(R/I^{[p^{n}]}) - \sum_{i = 1}^\ell \ell_{R}(R/I_i^{[p^{n}]})
\end{eqnarray*}
since ${}_{F^n}R$ is a Cohen-Macaulay $R$-module
and the residue class field of $R$ is algebraically closed.

We set $I_i = (a_{i1}, \ldots, a_{id})$, 
where $a_{i1}, \ldots, a_{id}$ forms a maximal $R$-regular sequence.
Then, 
\[
\ell_{R}(R/I_i^{[p^{n}]}) = \ell_{R}(R/(a_{i1}^{p^n}, \ldots, a_{id}^{p^n}))
= p^{dn} \ell_{R}(R/(a_{i1}, \ldots, a_{id})) .
\]

Thus, we have
\[
\ell_{R}(R/I^{[p^{n}]})
= \left( \epsilon_d \beta_d + 
\sum_{i = 1}^\ell \ell_{R}(R/(a_{i1}, \ldots, a_{id})) \right) p^{dn}
+ \sum_{i = 0}^{d-1} \epsilon_i \beta_i p^{in} .
\]

Remark that
\[
\epsilon_d \beta_d + 
\sum_{i = 1}^\ell \ell_{R}(R/(a_{i1}, \ldots, a_{id})) = e_{HK}(I) > 0.
\]
Putting $\alpha = e_{HK}(I)$, we know that $I$ satisfies the required condition. 
We have completed the proof of Theorem~\ref{HK}.

 In Theorem~\ref{HK}, the coefficients of the
polynomial are assumed to be zero in the terms of degree $d/2$ or
lower. This assumption is made due to the fact that
$\overline{{\rm A}_i(R)}_{\mathbb Q} = 0$ for $i\leq d/2$ 
(for a homogeneous coordinate ring $R$ of a smooth projective variety) if the
Grothendieck's standard conjecture holds ({\em c.f.}
\cite[Remark~7.12]{K23}). There is no known example where $\overline{{\rm A}_i(R)}_{\mathbb Q}$ does
not vanish for some  $i \leq d/2$. 

Therefore, it is natural to ask the following:

\begin{Conjecture}\label{yosou2}
Let $R$ be a $d$-dimensional Cohen-Macaulay local ring of characteristic $p$ with perfect residue class field.
Let $I$ be an maximal primary ideal of $R$ of finite projective dimension.
We set 
\[
\ell_R(R/I^{[p^n]}) = \sum_{i = 0}^d \beta_i p^{in} 
\]
for $n > 0$.
Then, if $i \le d/2$, $\beta_i = 0$.
\end{Conjecture}

\section{Proof of Lemma~\ref{lemma}}
This section is devoted to proving Lemma~\ref{lemma}.

We use the following basic properties on singular Riemann-Roch maps.

\begin{Fact}\label{3.1}
Let $X$ be a $d$-dimensional projective variety over $k$.
Then, we have an isomorphism
\[
G_0(X)_{\Bbb Q} \stackrel{\tau_X}{\longrightarrow} {\rm A}_*(X)_{\Bbb Q} .
\]
Let ${\mathcal M}$ be a coherent ${\mathcal O}_X$-module.
Put 
\[
\tau_X([{\mathcal M}]) = s_d + s_{d-1} + \cdots + s_0 ,
\]
where $s_i \in {\rm A}_i(X)_{\Bbb Q}$.

Let $D$ be a very ample divisor on $X$.
Put $S = k[x_0, x_1, \ldots, x_n]$, where $S$ is a graded polynomial ring with
$\deg(x_i) = 1$ for $i = 0, 1, \ldots, n$.
Let 
\[
X = \proj B \stackrel{i}{\hookrightarrow} {\Bbb P}^n = \proj{S}
\]
be the embedding corresponding to $D$,
where we have 
\[
S \twoheadrightarrow  B = k[B_1] \subset
\bigoplus_{m} H^0(X, {\mathcal O}_X(mD)) .
\]
Here, $B_1$ denotes the homogeneous component of the graded ring $B$ of degree $1$.
We note that $B$ is standard graded; that is, $B$ is a graded ring 
generated by elements of degree $1$ over $B_0 = k$.
\begin{enumerate}
\item
We have a commutative diagram:
\[
\begin{array}{ccc}
G_0(X)_{\Bbb Q} &  \stackrel{\tau_X}{\longrightarrow} & 
{\rm A}_*(X)_{\Bbb Q} \\
i_* \downarrow \phantom{i_*} & & \phantom{i_*} \downarrow i_* \\
G_0({\Bbb P}^n)_{\Bbb Q} &  \stackrel{\tau_{{\Bbb P}^n}}{\longrightarrow} & 
{\rm A}_*({\Bbb P}^n)_{\Bbb Q}
\end{array}
\]
Put
\[
M = \bigoplus_{m} H^0(X, {\mathcal M}\otimes_{{\mathcal O}_X}
{\mathcal O}_X(mD)) .
\]

Then, $M$ is a graded $\bigoplus_{m} H^0(X, {\mathcal O}_X(mD))$-module.
We have
\[
\tau_{{\Bbb P}^n}([\tilde{M}])  =  \tau_{{\Bbb P}^n}i_*([{\mathcal M}])  =  i_*\tau_X([{\mathcal M}])  =  i_*(s_d) + i_*(s_{d-1}) + \cdots + i_*(s_0) .
\]
Put
\[
H_i = \proj{S/(x_{i+1}, \ldots, x_n)} .
\]
Then, we have 
\[
{\rm A}_i({\Bbb P}^n)_{\Bbb Q} = {\Bbb Q}[H_i]
\]
for $i = 0, 1, \ldots, n$.
Let $\ell_i$ be a rational number such that
\[
i_*(s_i) = \ell_i[H_i]
\]
for $i = 0, 1, \ldots, d$.
Then, we have
\[
P_M(t) = \dim_kM_t = 
\frac{\ell_d}{d!}t^d + \frac{\ell_{d-1}}{(d-1)!}t^{d-1} +
\cdots + \frac{\ell_0}{0!}t^0 
\]
for $t \gg 0$.
(See Proposition in p3005 in Chan-Miller~\cite{CM}, Proposition~4.1 in Roberts-Singh~\cite{RSi})
\item
Let $m$ be the homogeneous maximal ideal of $B$.
We have the following commutative diagram:
\[
\begin{array}{ccc}
G_0(X)_{\Bbb Q} &  \stackrel{\tau_X}{\longrightarrow} & 
{\rm A}_*(X)_{\Bbb Q} \\
\alpha \downarrow \phantom{\alpha} & & \phantom{\beta} 
\downarrow \beta \\
G_0(B)_{\Bbb Q} &  \stackrel{\tau_{B}}{\longrightarrow} & 
{\rm A}_*(B)_{\Bbb Q} \\
\gamma \downarrow \phantom{\gamma} & & \phantom{\delta} 
\downarrow \delta \\
G_0(B_m)_{\Bbb Q} &  \stackrel{\tau_{B_m}}{\longrightarrow} & 
{\rm A}_*(B_m)_{\Bbb Q}
\end{array}
\]
The horizontal maps are isomorphisms.
Here, $\gamma$ and $\delta$ are isomorphisms induced by 
localization $B \rightarrow B_m$.
For a graded $B$-module $M$, we have
\[
\alpha([\tilde{M}]) = [M] .
\]
Here, we remark that $\alpha$ is well-defined
since $[T] = 0$ in $G_0(B)_{\Bbb Q}$ for a graded $B$-module $T$
whose homogeneous graded pieces are zero except for finitely many degrees. 
The map $\beta$ is the sum of the maps
\[
{\rm A}_i(X)_{\Bbb Q} \stackrel{\beta}{\twoheadrightarrow} 
\frac{{\rm A}_i(X)_{\Bbb Q}}{c_1(D){\rm A}_{i+1}(X)_{\Bbb Q}}
= {\rm A}_{i+1}(B)_{\Bbb Q} ,
\]
where this map is given by
\[
[\proj{B/P}] \mapsto [\spec(B/P)]
\]
for each homogeneous prime ideal $P$ with $\dim B/P > 0$.

Thus, we obtain
\[
\tau_{B_m}([M \otimes_BB_m])
= \delta\beta(s_d) + \delta\beta(s_{d-1}) + \cdots +
\delta\beta(s_0) ,
\]
where $\delta\beta(s_i) \in {\rm A}_{i+1}(B_m)_{\Bbb Q}$.
We refer the reader to section~4 in \cite{K11} for maps $\alpha$, $\beta$,
$\gamma$, $\delta$.
\end{enumerate}
\end{Fact}

\begin{Example}\label{symbolic}
Set $X = {\Bbb P}^m \times {\Bbb P}^n$.
Let $p_1$ and $p_2$ be the first and second projections, 
respectively.
Assume $m \ge n \ge 2$.

Then, we have
\[
G_0(X)_{\Bbb Q} \stackrel{\tau_X}{\longrightarrow}
{\rm A}_*(X)_{\Bbb Q} = {\Bbb Q}[a,b]/(a^{m+1}, b^{n+1}) ,
\]
where 
\[
{\rm A}_{m+n-c}(X)_{\Bbb Q} = \bigoplus_{i+j = c}{\Bbb Q}a^ib^j ,
\]
and $a = c_1(p_1^*{\mathcal O}_{{\Bbb P}^m}(1)) \in {\rm A}_{m+n-1}(X)$ and 
$b = c_1(p_2^*{\mathcal O}_{{\Bbb P}^n}(1)) \in {\rm A}_{m+n-1}(X)$.
We put
\[
{\mathcal O}_X(s,t) = p_1^*{\mathcal O}_{{\Bbb P}^m}(s)
\otimes_{{\mathcal O}_X} p_2^*{\mathcal O}_{{\Bbb P}^n}(t) .
\]

Put
\[
f(x) = \frac{x}{1-e^{-x}} .
\]

Then, we have
\begin{eqnarray*}
\tau_{X}([{\mathcal O}_X(s,t)]) & = & 
{\rm ch}(p_1^*{\mathcal O}_{{\Bbb P}^m}(s))
{\rm ch}(p_2^*{\mathcal O}_{{\Bbb P}^n}(t))
{\rm td}(\Omega_X^\vee) \\
& = & 
{\rm ch}(p_1^*{\mathcal O}_{{\Bbb P}^m}(s))
{\rm ch}(p_2^*{\mathcal O}_{{\Bbb P}^n}(t))
{\rm td}(p_1^*\Omega_{{\Bbb P}^m}^\vee) 
{\rm td}(p_2^*\Omega_{{\Bbb P}^n}^\vee) \\
& = & 
e^{sa}f(a)^{m+1}e^{tb}f(b)^{n+1} .
\end{eqnarray*}

Here, take a very ample divisor $a+b \in {\rm A}_{m+n-1}(X)$.
Then, the homogeneous coordinate ring $B$ is 
defined by all the $2$-minors of the generic $(m+1) \times (n+1)$-matrix.
In this case, the cycle in $X$ corresponding to $a^ib^j$ is the
closed subscheme defined by the ideal generated by the entries
in the top $i$ rows and the left $j$ columns.
Then, we have
\[
G_0(B)_{\Bbb Q} \stackrel{\tau_B}{\longrightarrow}
{\rm A}_*(B)_{\Bbb Q} = {\Bbb Q}[a,b]/(a^{m+1}, b^{n+1}, a+b)
= {\Bbb Q}[b]/(b^{n+1}) ,
\]
where we identify $a$ with $-b$.
Let $P$ (resp.\ $Q$) be the ideals generated by 
the elements in the first row (resp.\ the first colummn).
Then, for $s > 0$ and $t > 0$,
\begin{eqnarray*}
[P^{(s)}] & = & [P^s] = \alpha([{\mathcal O}_X(-s,0)]) \\
\mbox{$[Q^{(t)}]$} & = & [Q^t] = \alpha([{\mathcal O}_X(0,-t)]) .
\end{eqnarray*}
Here, for an ideal $I$, $I^{(s)}$ denotes the $s$-th symbolic power of $I$.
Then, 
\begin{eqnarray*}
\tau_B([P^{(s)}]) & = & e^{sb}f(-b)^{m+1}f(b)^{n+1} =
f(-b)^{m+1-s}f(b)^{n+1+s} \in {\Bbb Q}[b]/(b^{n+1})
\\
\tau_B([Q^{(t)}]) & = & f(-b)^{m+1}e^{-tb}f(b)^{n+1} =
f(-b)^{m+1+t}f(b)^{n+1-t} \in {\Bbb Q}[b]/(b^{n+1}) ,
\end{eqnarray*}
since $e^b = f(b)/f(-b)$.
Here,
\[
\{ P^{(m)}, P^{(m-1)}, \ldots, P^{(1)}, B, Q^{(1)}, \ldots, 
Q^{(n-1)}, Q^{(n)}  \}
\]
is the set of rank one maximal Cohen-Macaulay modules.
It is easily verified calculating local cohomology modules
of Segre products~\cite{GW}.
If there exists non-negative integers
$q_0$, $q_1$, \ldots, $q_{m+n}$ satisfying
\[
\sum_{k = 0}^{m+n}q_k > 0 \ \ \mbox{and}
\]
\[
\sum_{k = 0}^{m+n}q_kf(-b)^{1+k}f(b)^{m+n+1-k}
= (\sum_{k = 0}^{m+n}q_k) + b^{n+1}( \cdots ) ,
\]
then
\[
(P^{(m)})^{\oplus q_0} \oplus 
\cdots \oplus (P^{(1)})^{\oplus q_{m-1}} \oplus B^{\oplus q_{m}} \oplus (Q^{(1)})^{\oplus q_{m+1}} \oplus 
\cdots \oplus (Q^{(n)})^{\oplus q_{m+n}}
\]
is a $B$-test module.
If $q_m > 0$, then it contains $B$ as a direct summand.

The authors do not know whether a test module (having $B$ as a direct summand) exists or not in this case.
\end{Example}

In order to prove Lemma~\ref{lemma}, it is enough to show the following Claim.

\begin{Claim}\label{mainclaim}
Let $d$ be a positive integer and $p$ a prime number.
Let $k$ be an algebraically closed field of characteristic $p$.
If $d$ is even, then we put
\[
S = k[x_0, x_1, \ldots, x_{d/2}] \ \ \mbox{and} \ \ 
T = k[y_0, y_1, \ldots, y_{(d/2)-1}] .
\]
If $d$ is odd, then we put
\[
S = k[x_0, x_1, \ldots, x_{(d-1)/2}] \ \ \mbox{and} \ \ 
T = k[y_0, y_1, \ldots, y_{(d-1)/2}] .
\]
We think that $S$ and $T$ are graded rings with
${\rm deg}(x_i) = {\rm deg}(y_j) = 1$ for each $i$ and $j$.  Let
$\ell$ be a sufficiently large integer.  We denote by
$S \# T^{(\ell)}$ be the Segre product of $S$ and $T^{(\ell)}$, that
is, $S \# T^{(\ell)} = \oplus_{m \ge 0}(S \# T^{(\ell)})_m$ with
$(S \# T^{(\ell)})_m = S_m \otimes_k T_{m\ell}$ 
$($see \cite{GW}$)$.\footnote{ For a graded ring $T$, $T^{(\ell)}$ denotes the $\ell$-th
  Veronese subring of $T$.  Please do not confuse it with the symbolic
  power of ideals as in Example~\ref{symbolic}. }  Let $A$ be the
localization of $S \# T^{(\ell)}$ at the homogeneous maximal ideal.

Then, the ring $A$ satisfies the following conditions:
\begin{enumerate}
\item[$(1)$] The ring $A$ is a $d$-dimensional $F$-finite Cohen-Macaulay
  normal local domain and the residue class field of $A$ is
  algebraically closed.
\item[$(2)$]
${\rm A}_i(A)_{\Bbb Q} = 
\overline{{\rm A}_i(A)}_{\Bbb Q} = 
\left\{
\begin{array}{lll}
{\Bbb Q} & & (\frac{d}{2} < i \le d), \\
0 & & (\mbox{otherwise}).
\end{array}
\right.$
\item[$(3)$]
There exists a maximal Cohen-Macaulay $A$-module $M$ such that 
$\tau_A([A \oplus M]) \in {\rm A}_d(A)_{\Bbb Q}$.
\end{enumerate}
\end{Claim}

If $d$ is even, then we set $m = d/2$ and $n = d/2 - 1$.
If $d$ is odd, then we set $m = n = (d-1)/2$.
Let $\ell$ be a positive integer.
Then, $a + \ell b$ is a very ample divisor on 
$X = {\Bbb P}^m \times {\Bbb P}^n$.
Put $B = S \# T^{(\ell)}$.

Calculating local cohomologies of Segre products (see \cite{GW}),
(1) will be easily proved.

(2) will be proved by the method due to Roberts-Srinivas~\cite{RS}.
In fact, the rational equivalence on cycles on $X = {\Bbb P}^m \times {\Bbb P}^n$ coincides with the numerical equivalence.
Put $A = B_m$.
Then, by Theorem~7.7 in \cite{K23},
we have isomorphisms
\[
{\rm A}_*(B)_{\Bbb Q} \simeq {\rm A}_*(A)_{\Bbb Q}
 \simeq \overline{{\rm A}_*(A)}_{\Bbb Q} .
\]
We know 
\[
{\rm A}_i(B)_{\Bbb Q} = 
\left\{
\begin{array}{ll}
{\Bbb Q} & (d/2 < i \le d), \\
0 & (\mbox{otherwise})
\end{array}
\right.
\]
by \cite{K11}.

In the rest, 
using Theorem~\ref{main}, we shall prove (3).
We shall prove that 
\begin{equation}\label{aim}
\mbox{$p_i\overline{\tau_A}(C_{CM}(A)) = \overline{{\rm A}_i(A)}_{\mathbb R}$
 for $d/2 < i < d$.}
\end{equation}

Here, we define
\[
N_q = \bigoplus_{s \in {\Bbb Z}}H^0(X, {\mathcal O}_X(q+s, \ell s)) .
\]
and prove the following lemma:

\begin{Lemma}
For any $\ell > 0$, $N_q$ is a maximal Cohen-Macaulay $B$-module 
if $-m \le q \le 0$.
\end{Lemma}

\proof
We have
\[
N_q = S(q) \# T^{(\ell)} .
\]
Let $m_1$ (resp.\ $m_2$) be the homogeneous maximal ideal 
of $S$ (resp.\ $T^{(\ell)}$).

Then, $H^i_{m_1}(S(q))_s \neq 0$ if and only if
\[
\mbox{$i = 0$ and $s \ge -q$}
\]
or
\[
\mbox{$i = m+1$ and $s \le -q-m-1$} .
\]
Further, $H^i_{m_2}(T^{(\ell)})_s \neq 0$ if and only if
\[
\mbox{$i = 0$ and $s \ge 0$}
\]
or
\[
\mbox{$i = n+1$ and $s \le - \lceil \frac{n+1}{\ell} \rceil$} .
\]
Here $\lceil \frac{n+1}{\ell} \rceil$ denotes the minimal integer which is
bigger than or equal to $\frac{n+1}{\ell}$.
We refer the reader to \cite{GW} for local cohomologies of Segre products.
Therefore, $N_q$ is a maximal Cohen-Macaulay module if and only if
\[
\left\{
\begin{array}{l}
- \lceil \frac{n+1}{\ell} \rceil < -q \\
-q-m-1 < 0 .
\end{array}
\right.
\]
It is equivalent to 
\[
-m-1 < q < \lceil \frac{n+1}{\ell} \rceil .
\]
Therefore, if $-m \le q \le 0$, then $N_q$ is a maximal Cohen-Macaulay module.
\qed

We set
\[
h_{m,q}(x) = (x + q + m)(x + q + m-1) \cdots (x+q+1) .
\]

Consider the polynomials
\begin{eqnarray*}
h_{m,0}(x) & = & (x + m)(x + m-1) \cdots (x+2)(x+1), \\
h_{m,-1}(x) & = & (x + m-1)(x + m-2) \cdots (x+1)x, \\
h_{m,-2}(x) & = & (x + m-2)(x + m-3) \cdots x(x-1), \\
& \vdots &  \\
h_{m,q}(x) & = & (x+m+q)(x+m+q-1) \cdots (x+1+q), \\
& \vdots & \\
h_{m,-m}(x) & = & x(x-1)(x-2)\cdots (x - (m-1)) .
\end{eqnarray*}

The following lemma will be used later.

\begin{Lemma}\label{ind}
Suppose $m \ge 2$ and $m > u > 0$.
The set of the coefficients of $x^u$ in 
\[
h_{m,-1}(x), \ h_{m,-2}(x), \ \ldots, \ h_{m,-m}(x) 
\]
contains a negative value.
\end{Lemma}

\proof
We shall prove it by induction on $m$.

Suppose $m = 2$.
Then,
\begin{eqnarray*}
h_{2,-1}(x) & = & (x+1)x = x^2 +x , \\
h_{2,-2}(x) & = & x(x-1) = x^2 -x .
\end{eqnarray*}

Assume that $m \ge 2$ and the assertion is true for $m$.

Suppose $1 < u < m+1$.
By the induction hypothesis, there exists $-m \le q < 0$
such that the coefficient of $x^{u-1}$ in $h_{m,q}(x)$ is negative.
If the coefficient of $x^u$ in $h_{m,q}(x)$ is negative,
the coefficient of $x^u$ in
\[
h_{m+1,q}(x) = (x+q+m+1)h_{m,q}(x)
\]
is negative.
If the coefficient of $x^u$ in $h_{m,q}(x)$ is non-negative,
the coefficient of $x^u$ in
\[
h_{m+1,q-1}(x) = h_{m,q}(x)(x+q)
\]
is negative.

Suppose $u = 1$.
By the induction hypothesis, there exists $-m \le q < 0$
such that the coefficient of $x$ in $h_{m,q}(x)$ is negative.
Remark that $h_{m,q}(0) = 0$.
Then, the coefficient of $x$ in 
\[
h_{m+1,q}(x) = (x+q+m+1)h_{m,q}(x)
\]
is negative.
\qed

Consider
\[
\tau_X({\mathcal O}_X(q,0))
= e^{qa}f(a)^{m+1}f(b)^{n+1} \in {\rm A}_*(X)_{\Bbb Q} 
= {\Bbb Q}[a,b]/(a^{m+1}, b^{n+1}).
\]

\begin{Lemma}\label{lemma2}
Suppose that $v$ is an integer such that $1 \le v \le n$.
\begin{enumerate}
\item
Assume $v < m$.
Then, 
the set of the coefficients of $a^v$ in 
\[
\tau_X({\mathcal O}_X(-m,0)), \ 
\tau_X({\mathcal O}_X(-m+1,0)), \ \ldots, \ 
\tau_X({\mathcal O}_X(0,0))
\]
contains a positive value and a negative value.
\item
Assume $v = m = n$.
Then,  
the coefficient of $a^m$ in $\tau_X({\mathcal O}_X(0,0))$
is positive.
The coefficient of $a^m$ in $\tau_X({\mathcal O}_X(-1,0))$
is zero.
The coefficient of $a^{m-1}b$ in $\tau_X({\mathcal O}_X(-1,0))$
is positive.
\end{enumerate}
\end{Lemma}

\proof
The coefficient of $a^v$ in
\[
e^{qa}f(a)^{m+1}f(b)^{n+1}
\]
is equal to the coefficient of $a^v$ in
\[
e^{qa}f(a)^{m+1} .
\]
Since 
\[
\tau_{{\Bbb P}^m}([{\mathcal O}_{{\Bbb P}^m}(q)])
= e^{qa}f(a)^{m+1} \in {\Bbb Q}[a]/(a^{m+1}) ,
\]
the coefficient of $a^v$ in
\[
e^{qa}f(a)^{m+1} 
\]
is equal to
\begin{equation}\label{siki}
(m - v)! \left\{
\mbox{the coefficient of $x^{m-v}$ in the polynomial
$x + q + m \choose m$} 
\right\}
\end{equation}
by Fact~\ref{3.1} (1).
Furthermore, (\ref{siki}) is equal to
\[
\frac{(m - v)!}{m!} \left\{
\mbox{the coefficient of $x^{m-v}$ in the polynomial
$h_{m,q}(x)$}
\right\}  .
\]
It is easy to see that, for $0 < u < m$, the coefficient of $x^u$ in 
$h_{m,0}(x)$ is positive.
Therefore, Lemma~\ref{lemma2} (1) immediately follows from Lemma~\ref{ind}.

Assume that $v = m = n$.
Since the constant term of $h_{m,0}(x)$ is positive,
the coefficient of $a^m$ in
\[
f(a)^{m+1}f(b)^{m+1}
\]
is positive.
Since the constant term of $h_{m,-1}(x)$ is zero,
the coefficient of $a^m$ in
\[
e^{-a}f(a)^{m+1}f(b)^{m+1}
\]
is zero.
The coefficient of $a^{m-1}b$ in
\[
e^{-a}f(a)^{m+1}f(b)^{m+1}
\]
is equal to 
\[
\left\{ \mbox{the coefficient of $a^{m-1}$ in $e^{-a}f(a)^{m+1}$} 
\right\}
\times \left( \frac{m+1}{2} \right) ,
\]
where $\frac{m+1}{2}$ is 
the coefficient of $b$ in $f(b)^{m+1}$.
The sign of the coefficient of $a^{m-1}$ in 
$e^{-a}f(a)^{m+1}$ is the same as the sign of the coefficient 
of $x$ in $h_{m,-1}(x)$, that is obviously positive.
\qed

We return to the proof of Claim~\ref{mainclaim}\,(3) and take an ample
divisor $a + \ell b$ on $X = \mathbb P^m \times \mathbb P^n$ for $\ell > 0$.
Remark that $S \# T^{(\ell)}$ is a homogeneous coordinate ring of $X$ under the embedding 
corresponding to $a + \ell b$.
We denote this ring simply by $B$. 
Then the commutative diagram from Fact~\ref{3.1}(2) with the current  $X$ is precisely
\[
\begin{array}{cccl}
G_0(X)_{\Bbb Q} &  \stackrel{\tau_X}{\longrightarrow} & 
{\rm A}_*(X)_{\Bbb Q} & = {\Bbb Q}[a,b]/(a^{m+1},b^{n+1}) \\
\alpha \downarrow \phantom{\alpha} & & \phantom{\beta} 
\downarrow \beta & \\
G_0(B)_{\Bbb Q} &  \stackrel{\tau_{B}}{\longrightarrow} & 
{\rm A}_*(B)_{\Bbb Q} & = {\Bbb Q}[a,b]/(a^{m+1},b^{n+1},a+ \ell b) 
= {\Bbb Q}[b]/(b^{n+1})
\end{array}
\]
where $\beta(a) = -\ell b$.

Recall that 
\[
N_0, \ N_{-1}, \ldots, N_{-m}
\]
are graded Cohen-Macaulay $B$-modules.
Since $\overline{{\rm A}_i(A)}_{\mathbb R} = \mathbb R$, in order to show (\ref{aim}),  
it is enough to prove that, for $v = 1, 2, \ldots, n$,
the set of the coefficients of $b^v$ in 
\[
\tau_B([N_0]), \ \tau_B([N_{-1}]), \ldots, \tau_B([N_{-m}])
\]
contains a positive value and a negative value.
Note that 
\[
\tau_B([N_q]) = \tau_B\alpha ({\mathcal O}_X(q,0))
= \beta \tau_X ({\mathcal O}_X(q,0))
= \beta  \left( e^{qa}f(a)^{m+1} f(b)^{n+1}   \right) .
\]
Here recall that the map
\[
\beta : {\Bbb Q}[a,b]/(a^{m+1}, b^{n+1}) \longrightarrow
{\Bbb Q}[b]/(b^{n+1})
\]
is given by $\beta(a^sb^t) = (-1)^s \ell^s b^{s+t}$.
Thus, we have 
\[
\beta(\sum_{s, t \ge 0} q_{st}a^sb^t) = 
\sum_{s, t \ge 0} (-1)^sq_{st}\ell^sb^{s+t}
= \sum_{v = 0}^n(\sum_{s=0}^v(-1)^sq_{s,v-s}\ell^s)b^v .
\]
If $(-1)^vq_{v0} > 0$ (resp.\ $(-1)^vq_{v0} < 0$),
the coefficient of $b^v$ in the above is positive
(resp.\ negative) for $\ell \gg 0$.

First suppose $1 \le v \le n$ and $v < m$.
By Lemma~\ref{lemma2} (1), the set of coefficients of $b^v$ in 
\[
\beta  \left( f(a)^{m+1} f(b)^{n+1}   \right) , \ 
\beta  \left( e^{-a}f(a)^{m+1} f(b)^{n+1}   \right) , \ 
\ldots, \ 
\beta  \left( e^{-ma}f(a)^{m+1} f(b)^{n+1}   \right)
\]
contains a positive value and a negative value for $\ell \gg 0$.

Next suppose $v = m = n$.
By Lemma~\ref{lemma2} (2), the sign of 
the coefficients of $b^m$ in 
\[
\beta  \left( f(a)^{m+1} f(b)^{m+1}   \right)  \mbox{\ and \ }
\beta  \left( e^{-a}f(a)^{m+1} f(b)^{m+1}   \right) 
\]
are different for  $\ell \gg 0$.
We have complete the proof of Claim~\ref{mainclaim}.

\noindent
\begin{tabular}{l}
C-Y. Jean Chan \\
Department of Mathematics \\
Central Michigan University \\
Mt. Pleasant, MI 48858 \\
U. S. A. \\
{\tt chan1cj@cmich.edu} \\
{\tt http://people.cst.cmich.edu/chan1cj }
\end{tabular}

\vspace{2mm}

\noindent
\begin{tabular}{l}
Kazuhiko Kurano \\
Department of Mathematics \\
School of Science and Technology \\
Meiji University \\
Higashimita 1-1-1, Tama-ku \\
Kawasaki 214-8571, Japan \\
{\tt kurano@isc.meiji.ac.jp} \\
{\tt http://www.isc.meiji.ac.jp/\~{}kurano}
\end{tabular}

\end{document}